\documentclass[12pt]{article}%
\usepackage[applemac]{inputenc}
\usepackage{amsmath,amssymb,fullpage}
\usepackage{amsfonts}
\usepackage{amsmath,amssymb}
\usepackage[applemac]{inputenc}
\usepackage{amsmath,amssymb,fullpage}
\usepackage{color}
\usepackage{color, colortbl}
\usepackage{amsmath}
\usepackage{amssymb}
\usepackage{graphicx}
\usepackage{tikz}
\usepackage{geometry}
\usetikzlibrary{arrows}
\usepackage{amsfonts}
\usepackage{amsmath,amssymb}
\usepackage[applemac]{inputenc}
\usepackage{amsmath,amssymb,fullpage}
\usepackage{color}
\usepackage{amsmath}
\usepackage{amssymb}
\usepackage{graphicx}
\usepackage{tikz}
\newtheorem{theorem}{Theorem}[section]
\newtheorem{lemma}[theorem]{Lemma}
\newcommand{\E}{E}

\newtheorem{corollary}[theorem]{Corollary}
\newtheorem{definition}[theorem]{Definition}

\newtheorem{example}[theorem]{Example}

\newtheorem{problem}[theorem]{Problem}
\newtheorem{remark}[theorem]{Remark}

\usetikzlibrary{arrows,shapes,automata,petri}

\newcommand{\dproof}{\noindent {Proof.} \quad}
\newcommand{\fproof}{\hfill $\square$ \bigskip}

\numberwithin{equation}{section}

\def\RB{\mathbb{R}}

\definecolor{LightCyan}{rgb}{0.88,1,1}

\def\RR{{\mathbb{ R}}}

\def\EE{{\mathbb{ E}}}

\def\1B{\text{1\!\!I}}

\def\F{\mathcal{F}}
\def\E{\mathbb{E}}
\def\P{\mathbb{P}}
\def\R{\mathbb{R}}
\def\N{\mathbb{N}}
\def\d{\diamond}
\def\D{\mathcal{D}}
\def\P{\mathbb{P}}
\def\J{\mathcal{J}} 

\begin{document}

\title{Space-time stochastic calculus and white noise}
\author{Bernt \O ksendal$^{1}$}
\date{\textcolor{blue}{Dedicated to Catriona Byrne,\\
 in gratitude for her support and encouragement through 40 years }
 \vskip 0.5cm
17 May 2022}
\maketitle

\footnotetext[1]{%
Department of Mathematics, University of Oslo, Norway. 
Email: oksendal@math.uio.no.}
\begin{abstract}

In the first part of this paper I give the historical background to my initial interest in stochastic analysis and to the writing of my book Stochastic Differential Equations. The first edition of this book was published by Springer in 1985, with the highly appreciated support of Catriona Byrne.

In the second part I present a motivation for modelling the dynamics of a system subject to a noise by means of a stochastic partial differential equation (SPDE) driven by a time-space Brownian sheet.  This is followed by a brief survey of time-space white noise and Hida–Malliavin calculus, which are useful tools for studying such equations. 

As an illustration I apply white noise calculus to find an explicit solution of an SPDE describing population growth in an environment subject to time-space white noise.

\end{abstract}

%%%%%%%%%%%%%%%%%%%%%%%%%%%%%%%%%

\section{Introduction}

My interest in the interplay between mathematical analysis and probability theory goes back to the beginning of my studies at the University of Oslo in the mid 1960's. But it really gained momentum in the late 1970's when I started studying the beautiful little book ``Stochastic Integrals" by Henry McKean \cite{Mc}. My colleague at the time, Knut Aase, and I ran a little seminar on the book at Agder College (now the University of Agder)  and we were both fascinated by this new calculus that the book presented, namely the \emph{It\^o calculus! }\\

Then in 1982 Sandy Davie and Alan Sinclair offered me a Research Fellowship and invited me to spend one semester at the Department of Mathematics, University of Edinburgh. There they asked me to give a course on stochastic differential equations (SDEs). I knew nothing about SDEs, but started immediately to study the subject intensively, in an effort to at least stay ahead of my (very advanced) audience. It was a rewarding experience, which opened up a new world, consisting of both interesting new mathematics and a number of important applications, e.g. to modeling of dynamical systems with noise, filtering theory, optimal stopping, stochastic control and (subsequently) mathematical finance. I was enthusiastic about this new field of mathematics and started lecturing about it, almost like preaching a gospel, at a number of places, including Agder College, E\" otv\" os Lor\'and University, University of Troms\o, and California Institute of Technology (Caltech). Every time I used the opportunity to polish my lecture notes, based on useful feedback from the audience. \\

In 1984 I submitted my lecture notes to Springer, and asked them to consider the manuscript for publication, for example in the Springer Lecture Notes in Mathematics (LNM). At a conference in Lancaster the same year I was fortunate to meet the recently appointed Springer Mathematical Editor, namely \emph{Catriona Byrne} who was handling my manuscript. She saw the potential of the manuscript and recommend that it be considered for publication in the new book series Universitext. The reviews were all very positive, and in 1985 the first edition of my book ``Stochastic Differential Equations"  \cite{O1} appeared, and it became a great success.  I think one reason for the success of the book was that it filled a gap in the literature. There were already several excellent books available, written by top experts, but my book was written from the point of view of an enthusiastic beginner, who was not trying to humiliate the reader but to work with the reader to understand the topic. Some years later a professor once told me that although more recent and more polished editions of my book were available, he gave his students the first edition to start with, because it was somehow more ``raw'' and direct, without all the technical subtleties (that should be taken seriously later)  and therefore easier as a first encounter with the field.

I will always be grateful to Catriona for her encouragement and support, not just for the first edition of this book, but also for later editions \cite{O} and for the other books I wrote later with coauthors. I think she is an important reason for the success of Springer among the mathematical community.

\section{ Equations with noise and white noise theory}

One of the fascinations with stochastic analysis is the interplay between mathematical analysis and probability theory. And perhaps the most spectacular example of such interplay is the topic of dynamical systems subject to \emph{noise}. Noise in some form is everywhere in our society. For example, it can be in the form of mechanical noise from machines, noise due to lack of information in the system, or environmental noise due to random fluctuations in weather. A classic example is the model for population growth:\\
Let $Y(t)$ denote the density of a given population at time $t$. Then the most basic model for the growth of $Y(t)$ is the differential equation
\begin{align}
\frac{dY(t)}{dt}= \alpha_0(t) Y(t); \quad Y(0)=y_0 \text{ (constant)}\label{sde1}
\end{align}
where $\alpha_0(t)$ is a given function, representing the relative growth rate. 
The solution of this differential equation is
\begin{align*}
Y(t)=y_0 \exp(\int_0^t \alpha_0(s) ds).
\end{align*}

A weakness of this model is that it does not take into account that there might be unpredictable, random changes, or ``noise'' in the environment. If we add ``noise''  in the relative growth rate, the equation gets the form
\begin{align}
\frac{dY(t)}{dt}= \alpha(t) +\sigma(t) W(t); \quad Y(0)=y_0, \label{sde2}
\end{align}
where $W(t)$ represents ``noise''  at time $t$ and $\sigma$ is a given ``noise''  coefficient.
The question is what properties such a ``noise'' process $W(t)$ should have. If we try to represent $W(t)$ by a stochastic process defined on a filtered probability space $(\Omega, \F,\{\F_t\}_{t\geq 0},\P)$, then ideally, i.e. if we think of noise as ``white''  in some sense, we could require that 
$W(t_1)$ and $W(t_2)$ are independent if $t_1 \neq t_2$
 and that $W(t)$ is normalised, in the sense that
$\E[W(t)]=0$ (where $\E$ denotes expectation with respect to $\P$) and $\E[W^2(t)] < \infty$ for all $t$. But it turns out that no such measurable stochastic process exist. However, since the Brownian motion process $B(t)=B(t,\omega); t\geq 0, \omega \in \Omega$ is continuous and has stationary, independent increments, one could try to put
\begin{align}
W(t)=\frac{dB(t)}{dt}. \label{w1}
\end{align}
This derivative does not exist in the ordinary (strong) sense, but since $t \mapsto B(t,\omega)$ is continuous for a.a. $\omega$, it is weakly differentiable for a.a. $\omega$, and we could try to interpret it weakly, in the sense that we regard the equation \eqref{sde2}  as an integral equation, that is
\begin{align*}
Y(t)&=y_0 + \int_0^t Y(s)\Big(\alpha(s) +\sigma(s) \frac{dB(s)}{ds}\Big)ds \nonumber\\
&= y_0 + \int_0^t Y(s)\alpha(s)ds +\int_0^t \sigma(s) Y(s) dB(s).
\end{align*}
This last integral can be made rigorous as an \emph{It\^ o integral}. It is sometimes written in the following short hand, differential form
\begin{align*}
dY(t)= \alpha(t)Y(t) dt  +\sigma(t)Y(t) dB(t);\quad Y(0)=y_0.
\end{align*}
Applying the It\^ o rules of stochastic calculus (the It\^o formula), we find the following well-known solution
\begin{align}
Y(t)= y_0 \exp\Big(\int_0^t \sigma(s) dB(s) +\int_0^t \{\alpha(s)- \tfrac{1}{2}\sigma^2(s)\} ds\Big).\label{2.7}
\end{align}

\subsection{A more elaborate model}
In the model above we are only considering the population at one given point $x$ in space and the noise is in time only. In an attempt to get a more realistic model, let us consider the population density $Y(t,x)$ at the time $t \in \R$ and at the point $x\in \R^n$, where $n$ is the dimension of the space. For simplicity of the notation, let us assume that $n=1$ in the following.
In this case, it is natural to assume that the  ``noise''  depends on both time and space also, i.e. that it is represented by a 2-parameter process $W(t,x).$ This is also relevant in many other situations, for example in temperature modelling or, more generally, weather modelling.
Assuming this, and arguing as in \eqref{sde2} and allowing the coefficients to depend on both $t$ and $x$, we arrive at the following 2-parameter stochastic differential equation for $Y(t,x)$:
\begin{align*}
\frac{\partial^2Y(t,x)}{\partial t \partial x} =\alpha(t,x)Y(t,x)dt dx+\sigma(t,x) W(t,x)dt dx.
\end{align*}
Now suppose we proceed as above, and try to represent $W(t,x)$  weakly as
\begin{align}
W(t,x)= \frac{\partial^2 B(t,x)}{\partial t \partial x}, \label{w2}
\end{align}
where $B(t,x); t\geq0, x\in \R$ is a 2-parameter Brownian motion, also called a \emph{Brownian sheet}. (See the next section for explanation.)
Then we arrive at the following stochastic integral equation
\begin{align*}
Y(t,x)=y_0 +\int_0^t \int_0^x \alpha(s,z)Y(s,z) ds dz +\int_0^t \int_0^x \sigma(s,z) Y(s,z) B(ds,dz),
\end{align*}
where the last integral is the \emph{space-time (2-parameter) It\^o integral with respect to $B(\cdot,\cdot)$}, as constructed in \cite{C,W,Z}. See also \cite{N, P} and \cite{Wa}. As in the 1-parameter case, we will use the following short-hand differential notation
\begin{align}
\frac{\partial ^2 Y(t,x)}{\partial t \partial x}&= \alpha(t,x) Y(t,x)dt dx + \sigma(t,x) Y(t,x) B(dt,dx); \quad t,x \geq 0 \label{sde3}
\\ Y(t,0)&=Y(0,x)=y_0; \quad  \text{ for all } t,x \geq 0.\nonumber
\end{align}
It follows from Theorem 2.4.1 in \cite{N} that this stochastic partial differential equation (SPDE) has a unique adapted solution. The question is:

\begin{problem}

\emph{Can we find a formula for the solution of \eqref{sde3}, like we did in the 1-parameter case \eqref{2.7}? }
\end{problem}

In view of the many applications of white noise theory to 1-parameter stochastic calculus (see e.g.  \cite{AaOPU,AaOU,AO,AO2,DO,DOP,HOUZ,HO,HS,OR} and also the spectacular white noise solution of general SDEs in \cite{LaP}), it is natural to ask if such an interplay can be useful also in the multi-parameter case. To support this idea we will show that white noise theory  can be used to find the solution of \eqref{sde3}. We will explain how after a short review of \emph{multi-parameter white nose calculus}. It is mainly based on the presentation in \cite{HOUZ}, and we refer to that book for proofs and more details.

\section{Multi-parameter white noise calculus}
The basic idea of white noise analysis, due to Hida \cite{H}, is to
consider white noise $W$ rather than Brownian motion $B$ as the fundamental object.
Within this framework, we will see that Brownian motion can indeed be expressed as the integral of white
noise. 

\subsection{The  $d$-parameter white noise probability space}
 In the
following $d$ will denote a fixed positive integer, interpreted as
either the time-, space- or time-space dimension of the system we
consider. More generally, we will call $d$ the parameter
dimension. Let $\mathcal{S}(\R^d)$ be the Schwartz space of rapidly
decreasing smooth $(C^\infty)$ real valued functions on $\RR^d$. 
Recall that $\mathcal{S}(\R^d)$  is a Fr\'echet space under the
family of seminorms
$$\Vert
f\Vert_{k,\alpha}:=\sup\limits_{x\in\RR^d}\{(1+|x|^k)|\partial^\alpha
f(x)|\},$$
where $k$ is a nonnegative integer,
$\alpha=(\alpha_1,\cdots,\alpha_d)$ is a multi-index of
nonnegative integers $\alpha_1,\cdots,\alpha_d$ and
$$\partial^\alpha
f=\frac{\partial^{|\alpha|}f}{\partial x_1^{\alpha_1}\cdots\partial x_d^{\alpha_d}}
\quad{\rm
where}\quad|\alpha|:=\alpha_1+\cdots+\alpha_d.$$

The dual $\Omega={\cal S}'(\RR^d)$ of ${\cal S}(\RR^d)$, equipped
with the weak-star topology, is the space of {\it tempered
distributions}. This space will be the base of our basic
probability space, which we explain in the following:\\
As events we will use the family $\mathcal{F}={\cal
B}({\cal S}'(\RR^d))$ of Borel subsets of ${\cal S}'(\RR^d)$, and our
probability measure $\P$ is defined by the following result:

\begin{theorem}{\bf (The Bochner--Minlos theorem)}\\
There exists a unique probability measure $\P$ on ${\cal B}({\cal S}'(\RR^d))$
with the following property:
$$\E[e^{i\langle\cdot,\phi\rangle}]:=\int\limits_{\cal S'}e^{i\langle\omega,
\phi\rangle}d\mu(\omega)=e^{-\tfrac{1}{2} \Vert\phi\Vert^2};\quad i=\sqrt{-1}$$
for all $\phi\in{\cal S}(\RR^d)$, where
$\Vert\phi\Vert^2=\Vert\phi\Vert^2_{L^2(\RR^d)},\quad\langle\omega,\phi\rangle=
\omega(\phi)$ is the action of $\omega\in{\cal S}'(\RR^d)$ on
$\phi\in{\cal S}(\RR^d)$ and $\E=\E_{\P}$ denotes the expectation
with respect to  $\P$.
\end{theorem} 
We will call the triplet $({\cal S}'(\RR^d),{\cal
B}({\cal S}'(\RR^d)),\P)$ the {\it  white noise probability
space\/}, and $\P$ is called the {\it white noise probability measure}.
\vskip0.3cm
The measure $\P$ is also often called the (normalized) {\it Gaussian
measure\/} on ${\cal S}'(\RR^d)$. It is not difficult to prove that if $\phi\in L^2(\RR^d)$ and
we choose $\phi_n\in{\cal S}(\RR^d)$ such that $\phi_n\to\phi$ in $L^2(\RR^d)$,
then
$$\langle\omega,\phi\rangle:=\lim\limits_{n\to\infty}\langle\omega,\phi_n\rangle
\quad\text{exists in}\quad L^2(\P)$$
and is independent of the choice of $\{\phi_n\}$. In
particular, if we define
$$\widetilde{B}(x):=\widetilde{B}(x_1,\cdots,x_d,\omega)=\langle\omega,\chi_
{[0,x_1]\times\cdots\times[0,x_d]}\rangle; \quad  
x=(x_1,\cdots,x_d)\in\RR^d,$$
where $[0,x_i]$ is interpreted as $[x_i,0]$ if $x_i<0$,
then $\widetilde{B}(x,\omega)$ has an $x$-continuous version $B(x,\omega)$, which
becomes a \emph{$d$-parameter Brownian motion}, in the following sense:
\vskip0.3cm

By a {\it $d$-parameter Brownian motion\/} we mean a family
$\{B(x,\cdot)\}_{x\in\RR^d}$ of random variables on a probability space
$(\Omega,{\cal F},\P)$ such that
\begin{itemize} 
\item
$B(0,\cdot)=0\quad\text{almost surely with respect to } \P,$
\item
$\{B(x,\omega)\}$ is a continuous and Gaussian stochastic process 
\item
and, further, for all $x=(x_1,\cdots,x_d)$,
$y=(y_1,\cdots,y_d)\in\RR_+^d$,
$B(x,\cdot),\,B(y,\cdot)$ have the covariance
$\prod_{i=1}^dx_i\wedge y_i$. For general $x,y\in\RR^d$ the covariance is
$\prod_{i=1}^d\int_\RR\theta_{x_i}(s)\theta_{y_i}(s)ds$, where
$\theta_x(t_1,\dots,t_d)=\theta_{x_1}(t_1)\cdots\theta_{x_d}(t_d)$, with
\begin{equation*} 
\theta_{x_j}(s)=
\begin{cases}
1\quad \text{ if } 0<s\leq x_j\\
-1\quad \text{ if } x_j<s\leq 0\\
0 \quad \text{ otherwise}
\end{cases}
\end{equation*} 
\end{itemize} 
\

It can be proved that the process $\widetilde{B}(x,\omega)$ defined above has a modification $B(x,\omega)$ which satisfies all these properties.
This process $B(x,\omega)$ then becomes a \emph{$d$-parameter Brownian motion}. 
\vskip0.3cm
We remark that for $d=1$ we get the classical (1-parameter) Brownian motion
$B(t)$ if we restrict ourselves to $t\geq 0$. For $d \geq 2$ we get what is often
called the \emph{Brownian sheet}.
\vskip0.3cm
With this definition of Brownian motion it is natural to define the
$d$-parameter Wiener--It\^{o} integral of $\phi\in L^2(\RR^d)$ by
$$\int\limits_{\RR^d}\phi(x)dB(x,\omega):=\langle\omega,\phi\rangle;\quad
\omega\in{\cal S}'(\RR^d).$$
We see that by using the Bochner--Minlos theorem we have obtained an easy construction of
d-parameter Brownian motion that works for any parameter dimension $d$. Moreover, we get a representation of the space $\Omega$ as the Fr\' echet space $\mathcal{S}'(\R^d)$. This is an advantage in many situations, for example in the construction of the Hida-Malliavin derivative, which can be regarded as a stochastic gradient on $\Omega$. See Section 4 and e.g. \cite{DOP} and the references therein.

\subsection{Chaos expansion in terms of Hermite polynomials}

The {\it Hermite polynomials\/} $h_n(x)$ are defined by
$$h_n(x)=(-1)^n
e^{{1}/{2}x^2}\frac{d^n}{dx^n}(e^{-{1}/{2}x^2});\quad
n=0,1,2,\cdots.$$
We see that the first Hermite polynomials are
\begin{align*}
& h_0(x)=1,\;h_1(x)=x,\;h_2(x)=x^2-1,\;h_3(x)=x^3-3x,\\
& h_4(x)=x^4-6x^2+3,\;h_5(x)=x^5-10x^3+15x,\cdots.  
\end{align*}
The {\it Hermite functions} $\xi_n(x)$ are defined by
$$\xi_n(x)=\pi^{-{1}/{4}}((n-1)!)^{-{1}/{2}}e^{-{1}/{2}x^2}h_{n-1}
(\sqrt{2}x);\quad n=1,2,\cdots.$$
Some important properties of the Hermite functions are the following:
\begin{itemize}
\item
$\xi_n\in{\cal S}(\RR)\quad\text{for all}\;n$
\item
The collection $\{\xi_n\}^\infty_{n=1}$ constitutes an
orthonormal basis for $L^2(\RR).$
\item
$\sup\limits_{x\in\RR}|\xi_n(x)|=O(n^{-{1}/{12}}).$

\end{itemize}

\vskip0.3cm
We now use these functions to define an orthogonal basis for $L^2(\P)$:\\

In the following, we let
$\delta=(\delta_1,\cdots,\delta_d)$ denote $d$-dimensional multi-indices with
$\delta_1,\cdots,\delta_d\in\mathbb{N}$. It follows that the family of tensor
products
$$\xi_\delta:=\xi_{(\delta_1,\cdots,\delta_d)}:=\xi_{\delta_1}\otimes\cdots
\otimes
\xi_{\delta_d}\;;\;\delta\in\mathbb{N}^d$$
forms an orthogonal basis for $L^2(\RR^d)$. Let
$\delta^{(j)}=(\delta_1^{(j)},\delta_2^{(j)},\cdots,\delta_d^{(j)})$ be
the $j$th multi-index number in some fixed ordering of all $d$-dimensional
multi-indices $\delta=(\delta_1,\cdots,\delta_d)\in\N^d$. We may
assume that this ordering has the property that
$$i<j\Rightarrow\delta_1^{(i)}+\delta_2^{(i)}+\cdots+\delta_d^{(i)}\leq
\delta_1^{(j)}+\delta_2^{(j)}+\cdots+\delta_d^{(j)},$$
i.e., that the $\{\delta^{(j)}\}^\infty_{j=1}$ occur in increasing order.
\bigskip

Now define
$$\eta_j:=\xi_{\delta^{(j)}}=\xi_{\delta_1^{(j)}}\otimes\cdots\otimes\xi_{
\delta_d
^{(j)}};\;j=1,2,\cdots.$$
We regard multi-indices as elements of the space
$(\N_0^{\N})_c$ of all sequences $\alpha=(\alpha_1,\alpha_2,\cdots)$
with elements $\alpha_i\in\N_0$ and with compact support,
i.e., with only finitely many $\alpha_i\not=0$. Put
$${\cal J}=(\N_0^{\N})_c.$$

\begin{definition}
Let $\alpha=(\alpha_1,\alpha_2,\cdots)\in{\cal J}$. Then we define
$$H_\alpha(\omega):=\prod\limits_{i=1}^\infty
h_{\alpha_i}(\langle\omega,\eta_i\rangle);\quad\omega\in{\cal
S}'(\R^d).$$
\end{definition}

\begin{theorem}{(Wiener--It\^{o} chaos expansion theorem)} Every
$f\in L^2(\P)$ has a unique representation
$$f(\omega)=\sum\limits_{\alpha\in{\cal J}}c_\alpha
H_\alpha(\omega),$$
where $c_\alpha\in\R$ for all $\alpha$.
\end{theorem}

Moreover,  the following isometry holds:
$$\Vert
f\Vert^2_{L^2(\P)}=\sum\limits_{\alpha\in{\cal
J}}\alpha!c_\alpha^2.$$

\begin{example}
The $d$-parameter Brownian motion
$B(x,\omega)$  is defined by:
\vskip0.1cm
$$B(x,\omega)=\langle\omega,\psi\rangle,$$
where
$$\psi(y)=\chi_{[0,x_1]\times\cdots\times[0,x_d]}(y).$$
Proceeding as above, we write
$$\psi(y)=\sum\limits_{j=1}^\infty(\psi,\eta_j)\eta_j(y)=\sum\limits_{j=1}^
\infty \big(\int\limits_0^x\eta_j(u)du \big) \eta_j(y),$$
where we have used the multi-index notation
$$\int\limits_0^x\eta_j(u)du=\int\limits_0^{x_d}\cdots\int\limits_0^{x_1}
\eta_j(u_1,\cdots,u_d)du_1\cdots
du_d=\prod\limits_{k=1}^d\int\limits_0^{x_k}\xi_{\beta_k^{(j)}}(t_k)
dt_k$$

when $x=(x_1,\cdots,x_d)$. Therefore,
$$B(x,\omega)=\langle\omega,\sum\limits_{j=1}^\infty\int\limits_0^x\eta_j
(u)du\eta_j\rangle=\sum\limits_{j=1}^\infty \big(\int\limits_0^x\eta_j(u)du\big)
\langle\omega,\eta_j\rangle.$$
We conclude that $B(x,\omega)=B(x_1,x_2, ...,x_d,\omega)$ has the expansion
\begin{align}
B(x,\omega)&=\sum\limits_{j=1}^\infty\int\limits_0^x\eta_j(u)du\,H_ 
{\epsilon^{(j)}}(\omega)\nonumber\\
&=\sum\limits_{j=1}^\infty \Big(\int\limits_0^{x_1}\int\limits_0^{x_2} ... \int\limits_0^{x_d} \eta_j(u)du_1 du_2 ... du_d \Big)\,H_ 
{\epsilon^{(j)}}(\omega).\label{BM}
\end{align}
\end{example}

\subsection{The stochastic test function spaces $(\mathcal{S})$ and the stochastic distribution space $(\mathcal{S})^*$}

We have seen that the growth condition
\begin{align}
\sum\limits_\alpha\alpha!c_\alpha^2<\infty \label{3.1}
\end{align}
assures that
$$f(\omega):=\sum\limits_\alpha c_\alpha H_\alpha(\omega)\in L^2(\P).$$
By replacing condition \eqref{3.1}  by various other conditions
we will obtain a family of stochastic test functions and stochastic generalised functions that relates to
$L^2(\P)$ in a way that is analogous to the spaces ${\cal S}(\R^d)\subset
L^2(\R^d)\subset{\cal S}'(\R^d)$. These spaces provide a favorable setting for
the study of stochastic (ordinary and partial) differential equations. \\

As an analogue, recall the characterizations of ${\cal S}(\R^d)$ and ${\cal
S}'(\R^d)$ in terms of Fourier coefficients:
As above we let
$\{\delta^{(j)}\}^\infty_{j=1}=\{(\delta_1^{(j)},\cdots,
\delta_d^{(j)})\}^\infty_{j=1}$
be a fixed ordering of all $d$-dimensional multi-indices
$\delta=(\delta_1,\cdots,\delta_d)\in\N^d$. In general, if
$\alpha=(\alpha_1,\cdots,\alpha_j,\cdots)\in{\cal
J}$, $\beta=(\beta_1,\cdots,\beta_j,\cdots)\in(\R^\N)_c$ are two
finite sequences, we will use the notation
$$\alpha^\beta=\alpha_1^{\beta_1}\alpha_2^{\beta_2}\cdots\alpha_j^{\beta_j}
\cdots\quad\text{where}\quad\alpha_j^0=1.$$

\begin{theorem}{[Reed and Simon (1980), Theorem V. 13-14]}\\
{\bf a)} 
Let $\phi\in L^2(\R^d)$, so that
\begin{align}
\phi=\sum\limits_{j=1}^\infty a_j\eta_j, \label{phi}
\end{align}
where
$a_j=(\phi,\eta_j);\;j=1,2,\cdots,$
are the Fourier coefficients of $\phi$ with respect to 
$\{\eta_j\}^\infty_{j=1}$. Then
$\phi\in{\cal S}(\R^d)$ if and only if
\begin{align*}\sum\limits_{j=1}^\infty
a_j^2(\delta^{(j)})^\gamma<\infty,
\end{align*}
for {\it all\/} $d$-dimensional multi-indices
$\gamma=(\gamma_1,\cdots,\gamma_d)$.

{\bf b)} The space ${\cal S}'(\R^d)$ can be identified with the space of all
formal expansions
\begin{align}
T=\sum\limits_{j=1}^\infty b_j\eta_j \label{3.4}
\end{align}
such that
$$\sum\limits_{j=1}^\infty
b_j^2(\delta^{(j)})^{-\theta}<\infty$$
for {\it some\/} $d$-dimensional multi-index
$\theta=(\theta_1,\cdots,\theta_d)$.
\end{theorem}

If this condition holds, then the action of $T\in{\cal S}'(\R^d)$ given by
\eqref{3.4} on $\phi\in{\cal S}(\R^d)$ given by \eqref{phi} is
$$\langle T,\phi\rangle=\sum\limits_{j=1}^\infty a_j b_j.$$
\vskip0.3cm
We now formulate a stochastic analogue of this result. The following quantity
is crucial:
If $\gamma=(\gamma_1,\cdots,\gamma_j,\cdots)\in(\R^\N)_c$ (i.e., only
finitely many of the real numbers $\gamma_j$ are nonzero), we use the short-hand notation
$$(2\N)^\gamma:=\prod_{j=1}^{\infty} (2j)^{\gamma_j}.$$

\begin{definition}{\bf (The Hida spaces of stochastic test functions and stochastic
distributions)}

{\bf a)} {\bf The stochastic test function space $(\mathcal{S})$}\par
Let
$({\cal S})$
consist of the functions
$$f=\sum\limits_{\alpha \in \J} c_\alpha
H_\alpha\in L^2(\P)\quad\text{with}\;c_\alpha\in\R$$
such that
$$\Vert f\Vert^2_{k}:=\sum\limits_{\alpha \in \J}
c_\alpha^2(\alpha!)^{2}(2\N)^{k\alpha}<\infty\quad\text{for
all}\;k\in\N$$
equipped with the projective topology.\\

{\bf b)} {\bf The stochastic distribution spaces $(\mathcal{S})^*$}\par
Let $q\in \mathbb{R}$. We say that the formal sum $F=\sum\limits_{%
\alpha \in \mathcal{J}} b_{\alpha }H_{\alpha }$ belongs to the \emph{Hida
distribution Hilbert space} $(\mathcal{S})_{-q}$ if
\begin{equation}
\Vert F \Vert^2_{-q} :=\sum\limits_{\alpha \in \mathcal{J}}\alpha
!c_{\alpha}^{2} (2\mathbb{N})^{-\alpha q}<\infty .  \label{3.24}
\end{equation}
We define the \emph{Hida stochastic distribution space}%
\index{Hida distribution space} $(\mathcal{S})^{*}$ as the
union 
$(\mathcal{S})^{*} = \bigcup_{q\in\mathbb{R}} (\mathcal{S})_{-q}$
equipped with the inductive topology.
\end{definition}

Note that $(\mathcal{S})^{* }$ can be regarded as the dual of $(\mathcal{S})$%
as follows: \\
The action of $F=\sum\limits_{\alpha }b_{\alpha }H_{\alpha }\in (%
\mathcal{S})^{* }$ on $f=\sum\limits_{\alpha }a_{\alpha }H_{\alpha }\in (%
\mathcal{S})$, where $b_{\alpha}, a_{\alpha} \in \mathbb{R}$, is given by 
\begin{equation*}
\langle F,f\rangle =\sum\limits_{\alpha }\alpha !a_{\alpha }b_{\alpha }.
\end{equation*}
We have the inclusions 
\begin{equation*}
(\mathcal{S}) \subset L^{2}(\P) \subset (\mathcal{S})^{* }.
\end{equation*}

\begin{example}{\bf Singular white noise}\par
One of the most useful properties of $({\cal S})^\ast$ is that it contains the
{\it singular\/} or {\it pointwise white noise}:

The $d$-parameter singular white noise
process is defined by the formal expansion
\begin{align}W(x)=W(x,\omega)=\sum\limits_{k=1}^\infty\eta_k(x)H_{\epsilon^{(k)}}(\omega)
;\;x\in\R^d. \label{WN}
\end{align}

From the definition one van verify that
$$W(x,\omega)\in({\cal S})^{\ast}.$$

By comparing the expansion \eqref{WN} for singular white noise $W(x)$ and the expansion \eqref{BM} for Brownian motion $B(x)$, we see that
\begin{align}
W(x)={{\partial^d}\over{\partial x_1\cdots\partial x_d}}B(x)\hskip5mm\hbox{ in }({\cal S})^*. \label{WNBM}
\end{align}
In particular, for $d=1$, we put $x_1=t$ and get the familiar identity
$$W(t)={{d}\over{dt}}B(t)\hbox{ in }({\cal S})^*.$$
\end{example}

\subsection{The Wick product}

Since $x \mapsto B(x,\omega)$ is continuous a.s., it is weakly differentiable a.s., and we see that the identity \eqref{WNBM} also holds in the weak distribution sense (in $\mathcal{S}'(\R^d)$), a.s.  For that matter, one could argue that we might as well work on $\mathcal{S}'(\R^d)$ (a.s.) However, there is no tractable product operator on $\mathcal{S}'(\R^d)$. One of the (many) advantages of working on $(\mathcal{S})^{*}$ is that it has a natural multiplication  called the \emph{Wick product}, which is a binary operation on both $(\mathcal{S})$ and $(\mathcal{S})^{*}$, and is fundamental in the study of
stochastic (ordinary and partial) differential equations. In
general, one can say that the use of this product corresponds to
-- and extends naturally -- the use of It\^{o} integrals. We now
explain this in more detail.

\begin{definition}The {\it Wick product} $F\diamond G$ of two elements
$$F=\sum\limits_\alpha a_\alpha H_\alpha,\;G=\sum\limits_\alpha b_\alpha
H_\alpha\in({\cal S})^*\quad\text{with}\quad
a_\alpha,b_\alpha\in\R$$
is defined by
$$F\diamond
G=\sum\limits_{\alpha,\beta}(a_\alpha,b_\beta)H_{\alpha+\beta}.$$
\end{definition}

An important property of the spaces  $({\cal
S})^\ast,({\cal S})$ is that they are closed under Wick products:

\begin{lemma}
\parindent=1cm
\item{\bf a)} $F,G\in({\cal S})^{\ast}\Rightarrow F\diamond
G\in({\cal S})^{\ast}$;
\item{\bf b)} $f,g\in({\cal S})\Rightarrow f\diamond g\in({\cal
S})$.
\parindent=0cm
\end{lemma}

It is easy to see directly from the definition that the Wick product is commutative, associative  and distributive over addition.

The \it {Wick powers} $F^{\diamond k}\;;\;k=0,1,2,\cdots$ of $F\in({\cal
S})^{*}$ are defined inductively as follows:
\begin{align*}
&F^{\diamond 0}=1.\\
&F^{\diamond k}=F\diamond F^{\diamond(k-1)} \text{ for }
k=1,2,\cdots.
\end{align*}
The \emph{Wick exponential} of $F \in (\mathcal{S})^*$ is defined by
\begin{align*}
\exp^{\diamond} F = \sum_{n=0}^\infty \tfrac{1}{n!} F^{\diamond n}; \text{ if convergent in } (\mathcal{S})^*.
\end{align*}

\subsection{Wick product and Hermite polynomials}

There is a striking connection between Wick powers and Hermite polynomials $%
h_n; n=0,1,2, ...$:

\begin{theorem}
\label{th3.11} Choose $\varphi \in L^2(\mathbb{R}^d)$ and define the random variable 
$w_{\varphi}$ by 
$$w_{\varphi}(\omega)=\langle \omega, \varphi\rangle.$$
Then 
\begin{equation}
w_{\varphi} ^{\diamond n}= \|\varphi \|^{n} h_n(\tfrac{w_{\varphi}}{\| \varphi \|}) \label{herm} 
\end{equation}
where $\|\varphi\| = \|\varphi \|_{L^2(\mathbb{R}^d) }$.\newline

\end{theorem}

This result can be used to get an explicit formula for the Wick exponential:

\begin{theorem}{\bf The Wick exponential}

Let $\phi\in L^2(\R^d).$ Then
\begin{align*}
\exp^{\diamond}\Big(\int_{\RR^d} \phi(x) B(dx)\Big)=\exp\Big( \int_{\RR^d} \phi(x) B(dx)-\tfrac{1}{2} \int_{\RR^d} \phi^2(x)dx\Big).
\end{align*}
\end{theorem}

\dproof
We may assume that
$\phi=c\eta_1$, in which case we get, with $w(\phi)=\int_{\RR^d} \phi(x) B(dx)$,
\begin{align*}
\exp^\diamond[w(\phi)]&=\sum_{n=0}^\infty\frac{1}{n!}w(\phi)^{
\diamond
n}=\sum_{n=0}^\infty\frac{1}{n!}c^n\langle\omega,\eta_1\rangle^{\diamond n}\\
&=\sum_{n=0}^\infty\frac{c^n}{n!}H_{\epsilon_1}^{\diamond
n}(\omega)=\sum_{n=0}^\infty\frac{c^n}{n!}H_{n\epsilon_1}(\omega)\\
&=\sum_{n=0}^\infty\frac{c^n}{n!}h_n(\langle\omega,\eta_1\rangle)=\exp\Big[c
\langle\omega,
\eta_1\rangle-\tfrac{1}{2}c^2\Big]\\
&=\exp\bigg[w(\phi)-\tfrac{1}{2}\Vert\phi\Vert^2\bigg],
\end{align*}
where we have used the generating property of the Hermite polynomials.
\fproof
\subsection{Wick multiplication and It\^{o} integration}

One of the most
striking features of the Wick product is its relation to
It\^{o} integration. In short, this relation can be
expressed as follows:
\begin{theorem}
Let $Y(x)$ be an It\^ o integrable process. Then 
\begin{align}\int\limits_{\R^d}Y(x)
B(dx)=\int\limits_{\R^d}Y(x)\diamond W(x)dx. \label{wis}
\end{align}
\end{theorem}
Here the left hand side denotes the It\^o integral of the stochastic process
$Y(x)=Y(x_1,x_2, ... , x_d,\omega)$ with respect to $B(dx)=B(dx_1 dx_2 ... dx_d)$, while the right hand side is to be
interpreted as an $({\cal S})^\ast$-valued (Pettis) Lebesgue integral of $Y \d W$  in $(\mathcal{S})^*$.

\vskip0.3cm
This relation explains why the Wick product is so natural and important
in stochastic calculus. It is also the key to the fact that It\^{o} calculus (with
It\^{o}'s formula, etc.) with ordinary multiplication is equivalent to ordinary
calculus with Wick multiplication. To illustrate this,
consider the example with $d=1, x=t$ and $Y(t)=B(t)\cdot\chi_{_{[0,T]}}(t)$.
Then by the It\^{o} formula the left hand side of \eqref{wis} becomes

\begin{align}\int\limits_0^T
B(t)dB(t)=\frac{1}{2}B^2(T)-\frac{1}{2}T,\label{3.12} 
\end{align}
while (formal) Wick calculation makes the right hand side equal to
\begin{align*}
\int\limits_0^T B(t)\diamond W(t)dt=\int\limits_0^T B(t)\diamond
B'(t)dt=\frac{1}{2}B(T)^{\diamond 2},
\end{align*}
which is equal to \eqref{3.12}  in virtue of \eqref{herm}.

\section{The space-time Hida-Malliavin calculus}

It is natural to ask if also the Hida-Malliavin derivative can be extended to the space-time case. As in previous sections we assume that the Brownian motion $B(x)$; $x\in 
\mathbb{R}^d$, is constructed on the space $(\Omega ,\mathcal{B},\P)$ with $%
\Omega =\mathcal{S}^{\prime }(\mathbb{R}^d)$. Note that any $\gamma \in L^{2}(%
\mathbb{R}^d)$ can be regarded as an element of $\Omega =\mathcal{S}^{\prime }(%
\mathbb{R}^d)$ by the action 
\begin{equation*}
\left\langle \gamma ,\phi \right\rangle =\int_{\mathbb{R}^d}\gamma (x)\phi
(x)dx;\quad \phi \in \mathcal{S}(\mathbb{R}^d).
\end{equation*}

Following the approach in \cite{DOP} we define the Hida-Malliavin derivative as follows:

\begin{definition}
\label{direct.derivative}

\begin{description}
\item[(i)] Let $F\in L^{2}(\P)$ and let $\gamma \in L^2(\mathbb{R}^d)$ be
deterministic. Then the \emph{directional derivative of $F$} in $(\mathcal{S}%
)^*$ in the direction $\gamma$ is defined by 
\begin{equation}  \label{6.1}
\D_\gamma F(\omega) = \lim_{\varepsilon \to 0} \frac{1}{\varepsilon} \big[ %
F(\omega + \varepsilon \gamma) - F(\omega)\big]
\end{equation}
whenever the limit exists in $(\mathcal{S})^*$.

\item[(ii)] Suppose there exists a function $\psi: \mathbb{R}^d \mapsto (%
\mathcal{S})^*$ such
that 
\begin{equation}  \label{6.2}
\begin{split}
& \int_{{\mathbb{R}^d}} \psi(x) \gamma(x)dx \quad \text{ exists in } (%
\mathcal{S})^{*} \text{ and} \\
&\D_\gamma F = \int_{\mathbb{R}^d} \psi(x) \gamma(x) dx, \quad \text{ for all }
\gamma \in L^2(\mathbb{R}^d).
\end{split}%
\end{equation}
Then we say that $F$ is \emph{Hida-Malliavin differentiable} in $(\mathcal{S}%
)^*$ and we write 
\begin{equation*}
\psi (x) = D_x F, \quad x \in \mathbb{R}^d.
\end{equation*}
We call $D_x F \in (\mathcal{S})^*$  the \emph{Hida-Malliavin derivative of $F$ at $x$}.
\end{description}
\end{definition}

\begin{example}
Suppose $F(\omega )=\left\langle \omega ,f\right\rangle =\int_{%
\mathbb{R}^d}f(x)B(dx)$, $f\in L^{2}(\mathbb{R}^d)$. Then 
\begin{equation*}
\D_{\gamma }F=\frac{1}{\varepsilon }\big[\left\langle \omega +\varepsilon
\gamma ,f\right\rangle -\left\langle \omega ,f\right\rangle \big]%
=\left\langle \gamma ,f\right\rangle =\int_{\mathbb{R}^d}f(x)\gamma (x)dx.
\end{equation*}%
Therefore $F$ is Hida-Malliavin differentiable and 
\begin{equation*}
D_{x}\Big(\int_{\mathbb{R}^d}f(u)B(du)\Big)=f(x) \text{ for a.a. } x \in \R^d.
\end{equation*}
\end{example}
As in the 1-parameter case one can prove the following:

\begin{lemma}{(Chain rule)}\\
Let $F\in L^{2}(\P)$ be Hida-Malliavin differentiable, with $D_x F \in 
L^{2}(\lambda \times \P) $, where $\lambda$ denotes Lebesgue measure. Suppose that $\varphi \in C^1(\mathbb{R})$ and $%
\varphi^{\prime }(F) D_x F \in L^{2}( \lambda \times \P)$. Then $\varphi(F)$ is also Hida-Malliavin differentiable and 
\begin{equation}  \label{chain.rule}
D_x\big(\varphi(F)\big) = \varphi^{\prime }(F) D_x F\quad \text{ for a.a. } x \in \R^d.
\end{equation}
\end{lemma}

\subsection{The general Hida-Malliavin derivative}

The Hida-Malliavin derivative can be expressed in
terms of the Wiener-It\^o chaos expansion as follows:
 Recall that 
$$H_\alpha(\omega):=\prod\limits_{i=1}^\infty
h_{\alpha_i}(\theta_i);\quad \text{ where } \theta_i=\langle\omega,\eta_i\rangle$$
Then by the chain rule \eqref{chain.rule} we have, 
\begin{equation}  \label{6.6}
D_x H_\alpha = \sum_{k=1}^m \prod_{j\ne k} h_{\alpha_j}(\theta_j) \alpha_k
h_{\alpha_k -1}(\theta_k) \eta_k (x) = \sum_{k=1}^m \alpha_k \eta_k(x) H_{\alpha -
\epsilon^{(k)}}.
\end{equation}
In view of this, the following definition is natural:

\begin{definition}
\textbf{The general Hida-Malliavin derivative.} \label{HM}\newline
If $F=\sum_{\alpha \in \mathcal{J}}c_{\alpha }H_{\alpha }\in (\mathcal{S}%
)^{\ast }$ we define the \emph{Hida-Malliavin derivative} $D_{x}F$ of $F$ at 
$x$ in $(\mathcal{S})^{\ast }$ by the following expansion: 
\begin{equation}
D_{x}F=\sum_{\alpha \in \mathcal{J}}\sum_{k=1}^{\infty }c_{\alpha }\alpha
_{k}\eta_{k}(x)H_{\alpha -\epsilon ^{(k)}},  \label{hida7.7}
\end{equation}%
whenever this sum converges in $(\mathcal{S})^{\ast }$. 
\end{definition}

This extension of the Hida-Malliavin derivative makes it possible to deal with more general cases. 
In short, taking the Hida-Malliavin derivate $D_x F$ may take you from $L^2(\P)$  into $(\mathcal{S})^{*}$, but conditioning with respect to $\mathcal{F}_x$ brings you back to $L^2(\lambda \times \P)$. For example, in the case $d=1$ the following extension of the Clark-Ocone representation theorem \cite{O}  was proved in  \cite{AaOPU}:

\begin{theorem}
\textbf{(Generalised Clark-Ocone theorem ($d$=1))} \\
Let $F\in L^{2}(\mathcal{F}_{T},\P)$%
. Then $D_{t}F\in (\mathcal{S})^{\ast }$ for all $t$, $\mathbb{E}[D_{t}F|%
\mathcal{F}_{t}]\in L^{2}(\lambda \times \P)$ and 
\begin{equation*}
F=\mathbb{E}[F]+\int_{0}^{T}\mathbb{E}[D_{t}F|\mathcal{F}_{t}]dB(t).
\end{equation*}
\end{theorem}

It is natural to ask if a similar result can be obtained in the general parameter case:
\begin{problem}
Is there a $d$-parameter version of Theorem 4.5?
\end{problem}

\section{Solving the population growth equation}
We now have the machinery from white noise theory we need to solve the space-time stochastic partial differential equation \eqref{sde3}:\\
\begin{theorem}
Suppose $\alpha(t,x)$ and $\sigma(t,x)$ are deterministic and in $L^2(\R^2)$. Then the  solution $Y(t,x)$ of the equation
\begin{align*}
&\frac{\partial ^2 Y(t,x)}{\partial t \partial x}= \alpha(t,x) Y(t,x)dt dx + \sigma(t,x) Y(t,x) B(dt,dx); \quad t,x \geq 0,\nonumber\\
&Y(t,0)=Y(0,x)=y_0 \text{ (constant) for all }t,x\geq 0,
\end{align*}
can be written as
\begin{align*}
Y(t,x)=y_0 \sum_{n=0}^\infty \frac{1}{n!} \sum_{k=0}^n \tfrac{\| \sigma\|^{k}}{k! (n-k)!}\Big[ \Big(\int_0^t \int_0^x \alpha(s,z) ds dz\Big)^{n-k}  h_k \Big(\int_0^t \int_0^x \tfrac{\sigma(s,z)}{\|\sigma\|} B(ds,dz)\Big)\Big].
\end{align*}

\end{theorem}

\dproof
First note that the equation can be written
\begin{align}
Y(t,x)=y_0 + \int_0^t \int_0^x K(s,z) \diamond Y(s,z) ds dz, \label{4.1}
\end{align}
where
\begin{align}
K(s,z)=\alpha(s,z) + \sigma(s,z) W(t,x), \text{ with } W(t,x)=\frac{\partial^2}{\partial t \partial x} B(t,x).\label{4.2}
\end{align}
Substituting for $Y(s,z)$ in \eqref{4.1} we get
\begin{align*}
Y(u)= y_0 + \int_0^{u} K(u_1) \d \big(y_0+ \int_0^{u_1}K(u_2)\d Y(u_2)du_2\big)du_1,
\end{align*}
where we put $u=(t,x), u_j=(s_j,z_j); j=1,2, ...$.\\
Repeating this we obtain
\begin{align}
&Y(u)=y_0\Big( 1 + \int_0^{u} K(u_1)du_1 + \int_0^{u} (\int_0^{u_1} K(u_1) \d K(u_2)du_2 ) du_1 \nonumber\\
& +  ... + \int_0^{u} ( \int_0^{u_1} ... \int_0^{u_n} K(u_1)\d K(u_2) \d ... K(u_{n+1}) du_1 du_2 ...du_n) du_{n+1}\Big) +R_{n} \label{4.5}
\end{align}
where
\begin{align*}
R_{n}= y_0\int_0^{u} \big( \int_0^{u_1} ... \int_0^{u_n} K(u_1)\d K(u_2) \d ... K(u_{n+1}) \d Y(u_{n+1})du_1 du_2 ...du_n \big) du_{n+1}
\end{align*}
It follows from the V\aa ge inequality \cite{HOUZ},Theorem 3.3.1, that $R_{n} \rightarrow 0$ in $(\mathcal{S})^*$ as $n \rightarrow \infty$.
Therefore we get, by letting $n \rightarrow \infty$ in \eqref{4.5},
\begin{align*}
&Y(t,x)= y_0 \sum_{n=0}^\infty \int_0^{u} \big( \int_0^{u_1} ... \int_0^{u_n} K(u_1)\d K(u_2) \d ... K(u_{n+1}) du_1 du_2 ...du_n \big) du_{n+1}\nonumber\\
&=y_0 \sum_{n=0}^\infty \frac{1}{n! n!} \int_{[0,u]^n} K(u_1)\d K(u_2) \d ... \d K(u_n) (v)du_1 ... du_n\nonumber\\
&=y_0 \sum_{n=0}^\infty \frac{1}{n!n!} \Big(\int_{[0,u]} K(v) dv \Big)^{\d n}\nonumber\\
&=y_0 \sum_{n=0}^\infty \frac{1}{n!n!} \sum_{k=0}^n \tfrac{n!}{k! (n-k)!}\Big(\int_0^t \int_0^x \alpha(s,z) ds dz\Big)^{\d (n-k)} \Big(\int_0^t \int_0^x \sigma(s,x)B(ds,dz)\Big)^{\d k}\nonumber\\
&=y_0 \sum_{n=0}^\infty \frac{1}{n!n!} \sum_{k=0}^n \tfrac{n!}{k! (n-k)!} \Big(\int_0^t \int_0^x \alpha(s,z) ds dz\Big)^{n-k} \|\sigma\|^{k} h_k \Big(\int_0^t \int_0^x \tfrac{\sigma(s,z)}{\|\sigma\|} B(ds,dz)\Big)\nonumber\\
&=y_0 \sum_{n=0}^\infty \frac{1}{n!} \sum_{k=0}^n \tfrac{\| \sigma\|^{k}}{k! (n-k)!} \Big(\int_0^t \int_0^x \alpha(s,z) ds dz\Big)^{n-k}  h_k \Big(\int_0^t \int_0^x \tfrac{\sigma(s,z)}{\|\sigma\|} B(ds,dz)\Big).
\end{align*}
\fproof

\begin{remark}
A direct consequence of this result is that
\begin{align*}
\E [Y(t,x)]= y_0 \sum_{n=0}^\infty \frac{\|\sigma\|^{n}}{n!n!} \Big(\int_0^t \int_0^x \alpha(s,z) ds dz\Big)^{n}, 
\end{align*}
and if we assume that $\alpha=0$ , we get
\begin{align*}
Y(t,x)&=y_0 \sum_{n=0}^\infty \frac{\|\sigma\|^{n}}{n!n!} h_n\Big(\int_0^t \int_0^x \tfrac{\sigma(s,z)}{\|\sigma\|} B(ds,dz)\Big).
\end{align*}
\end{remark} 
\begin{remark}
It follows from the proof that this result holds for any random $y_0 \in L^2(\P)$, if we replace the product by a Wick product and interpret the equation in the Wick-It\^o-Skorohod sense, that is
\begin{align*}
&\frac{\partial ^2 Y(t,x)}{\partial t \partial x}= \alpha(t,x) Y(t,x)dt dx + \sigma(t,x) Y(t,x) \d W(t,x) dtdx; \quad t,x \geq 0,\nonumber\\
&Y(t,0)=Y(0,x)=y_0 \in L^2(\P) \text{ for all }t,x\geq 0.
\end{align*}
Then the solution is
\begin{align*}
Y(t,x)=y_0 \d \sum_{n=0}^\infty \frac{1}{n!} \sum_{k=0}^n \tfrac{\| \sigma\|^{k}}{k! (n-k)!} \Big(\int_0^t \int_0^x \alpha(s,z) ds dz\Big)^{n-k}  h_k \Big(\int_0^t \int_0^x \tfrac{\sigma(s,z)}{\|\sigma\|} B(ds,dz)\Big), 
\end{align*}
which has expectation
\begin{align*}
\E [Y(t,x)]= \E[y_0 ] \sum_{n=0}^\infty \frac{\|\sigma\|^{n}}{n!n!} \Big(\int_0^t \int_0^x \alpha(s,z) ds dz\Big)^{n}.
 \end{align*}
Note that no adaptedness conditions are needed.
\end{remark}

\section{Concluding remarks}
The multi-parameter stochastic calculus is a mostly unexplored area of research. It is clearly crucial in the study of stochastic partial differential equations driven by space-time Brownian motion. See e.g. \cite{Wa, N, NR, NZ} and also \cite{HOUZ}, which includes SPDEs driven by space-time Poisson random measure white noise as well. But it is also important for many other applications.  In this informal note I have tried to illustrate that white noise calculus can be a powerful tool in the study of multi-parameter stochastic calculus in general.

\section{Acknowledgments}
I am grateful to Nacira Agram and Yaozhong Hu for valuable comments.

\end{document}